\documentclass[a4paper,12pt]{article}
\begin{document}
\begin{center}
\bf A Pierce decomposition for\\
 generalized Jordan
triple systems of second order
\end{center}
\vskip 3mm
\begin{center}
Issai Kantor\\
Department of Mathematics, Lund University, S-221 00 Sweden\\
and\\
Noriaki Kamiya\\
Center for Mathematical
Sciences, University of Aizu, 965-8580 Japan\\
\end{center}
\vskip 3mm
\begin{center}
\bf Introduction
\end{center}
\par
In this paper we study a Peirce decomposition for
generalized Jordan
triple systems of second order.
\par
Let $(x,y,z)$ be a product in a 
generalized Jordan
triple system of second order
([11],[12])
and an element 
$e$ be a tripotent,\ i.e.,
$$
(eee)=e.$$
\par
Denote 
$$
L(x)=(eex),\ R(x)=(xee),\ Q(x)=(exe).$$
We prove that the space $U$ of a 
generalized Jordan
triple system of second order
is decomposed into a 
direct sum of eight components, such that each component consists of
eigen vectors of the linear operators $L,R$ whereas the
 action  of the linear operator $Q$  is somewhat more complicated. Namely,
$$
U=U_{00}\oplus U_{{1\over 2}{1\over 2}}\oplus
U_{11}\oplus
U_{{3\over 2}{3\over 2}}\oplus
U_{-{1\over 2}0}\oplus
U_{01}\oplus U_{{1\over 2}2}\oplus U_{13},\eqno(0.1)$$
when
$$
L(a)=\lambda a,\ R(a)=\mu a\ {\rm if}\ a\in U_{\lambda\mu}.
\eqno(0.2)
$$  
The operator  $Q$ acts as follows:  subspaces  $U_{11}$ and $U_{13}$
are presented as $U_{1 \mu}=U_{1 \mu}^{+}\oplus
U_{1 \mu}^{-}$ and 
$$
Q(a)=\pm  \mu a\ {\rm if}\ a\in U_{1 \mu}^\pm.
\eqno(0.3)
$$ 
Moreover 
$Q$ defines a one-to-one correspondence between the
subspaces
$U_{{3\over 2}{3\over 2}}
$ and $U_{{1\over 2}2}
$
(see 
Theorem 1.3) and $Q=0$ on the remaining subspaces.

This decomposition generalizes the Peirce decomposition for a Jordan
triple system [14], which consists only of the three first
components and could be written as
$$
U=U_{0}\oplus U_{1\over 2}\oplus U_{1},\eqno(0.3)
$$ 
because 
$R=L$ in the case of a Jordan triple system.
\par
We also consider a special case of a weakly commutative 
generalized Jordan
triple system of second order ([13]).
In that case the decomposition consists only of 
six components
(Theorem 1.3):
$$
U=U_{00}\oplus U_{{1\over 2}{1\over 2}}\oplus
U_{11}^{+}\oplus U_{11}^{-}\oplus
U_{01}\oplus 
U_{13}^{-}.\eqno(0.5)$$
\par
Another special case was in fact considered by   B. N. Allison in the paper [2].
  B. N. Allison gave the Peirce decomposition for the 
structurable
algebra defined by a hermitian idempotent satisfying one additional relation.
It is easy to see
  that such idempotent is also a tripotent of the 
corresponding triple system and the  Peirce decomposition for the algebra is the  Peirce decomposition for the  triple system in the case of these 
special tripotents.
In this case the decomposition has five  components :
$$U=U_{00}\oplus U_{{1\over 2}{1\over 2}}\oplus
U_{11}^{+} \oplus U_{11}^{-} \oplus  
U_{01}
.\eqno(0.6)$$
\par
The proof of Theorem 1.3 is based on the solution of a system of
all linear consequences of two identities in the 
definition of Jordan triple system of second order.
The solution of this system is given in Theorems 1.1 and
1.2.
The Theorems 1.1-1.3 are  the  content of \S 1.
\par
In \S 2 we give several examples of the Peirce decomposition.
\par
In \S 3 we go further and consider the
system of all bilinear consequences of identities of the
generalized Jordan triple system of second order.
We solve this system in the case when the tripotent is a left
unit.
The solution gives a set of relations
between components of the Peirce decomposition,
which we summarize in  Theorem 3.3.
\par
In the second author's paper [10],
a certain Peirce decomposition was investigated for the 
space of a Jordan triple system associated with the given
triple system $U$ but not a decomposition of the space of 
$U$ itself.
\par
Throughout this article, we consider triple systems over
a field $\Phi$
of characteristic
$\not= 2,3,5.$
\par
\vskip 3mm
\begin{center}
\bf \S 1\ A Peirce decomposition 
defined by a tripotent
\end{center}
\par
\vskip 3mm
{\bf Definition 1.}\quad
A vector space $U$
over a field $\Phi$
equipped with a triple product
$(xyz)$ is called a {\it generalized Jordan
triple system of second order}
 if
$$
(ab(cdf))=((abc)df)-(c(bad)f)+(cd(abf))\eqno(1.1)$$
$$
[((avb)uc)-(cu(avb))
-(cv(aub))-(a(vcu)b)]_{a,b}=0
\eqno(1.2)$$
where $[\;\; ]_{a,b}$  means alternation on
$a$ and $b.$

{\bf Remark.} The notion of generalized Jordan triple system (of the
finite order, particularly of the second order) was
introduced in the papers ([11], [12]) by the first author. Some authors are using
the terminology  {\it Kantor triple system} or  {\it Kantor pair} in
the case of the second order system ([3],[5]).
\par
Let $e$ be a tripotent, i.e.
$$
(eee)=e.$$
\par
Denote
$$
L(x)=(eex),\;Q(x)=(exe),\;R(x)=(xee).\eqno(1.3)$$
\par
We will find all conditions on these three operators
$L,Q,R$ derived from identities
(1.1) and (1.2).
For this purpose
we substitute in (1.1) and (1.2) instead of five letters $a,b,c,d,f$
four times $e$
and one time $x$.
The first four identities below realize all
such possibilities for the identity (1.1).
The equation
(1.8) is obtained by substitution in (1.2)
$a=x$ and  $e$ for other letters.
\par
We have
$$
R^{2}-Q^{2}+LR-R=0,\eqno(1.4)$$
$$
RQ-QR+LQ-Q=0,\eqno(1.5)$$
$$
LR=RL,\eqno(1.6)$$
$$
LQ+QL-2Q=0,\eqno(1.7)$$
$$
(R-2L-1)(R-L)=0.\eqno(1.8)$$

Moreover (1.6) and (1.8) imply
$$
(R-L)(R-2L-1)=0.\eqno(1.8^{'})$$

Our goal is to solve the system of
equations (1.4)-(1.8)
with three unknown linear operators
$L,R,Q.$
\par
\vskip 3mm
{\bf Lemma 1.1.}\quad
{\it There is no vector $a\not=0$
satisfying}
$$
La=Ra=-a.
\eqno(1.9)$$
\par
\vskip 3mm
{\it Proof.}\quad
If so, then the equations 
(1.4), (1.5) and (1.7) give
$$
Q^{2}a=3a,\ (R+L)Qa=0,\ LQa=3Qa.
\eqno(1.10)$$
\par
It follows from the first equality, 
  that
$b=Qa\not=0.$
The other two equalities give
$$
Lb=-Rb=3b.\eqno(1.11)$$
\par
Acting by (1.8) on $b$ and using (1.11),
one comes to a contradiction
$$
(-10)(-6)b=0.$$
\par
The lemma is proved.
\par
We denote
\begin{eqnarray*}
U_{R=L}&=&\{x\in U\;|\;Rx=Lx\},\\
U_{R=2L+1}&=&\{x\in U\;|\;Rx=(2L+1)x\},
\end{eqnarray*}
where
$U$ is the space of the triple system.
\par
\vskip 3mm
{\bf Lemma 1.2.}\quad
{\it The space $U$ is a direct sum of  subspaces
$U_{R=L}$ and $U_{R=2L+1}:$}
$$
U=U_{R=L}\oplus U_{R=2L+1}.
\eqno(1.12)$$
\par
\vskip 3mm
{\it Proof.}\quad
It follows from (1.8) and $(1.8^{'})$ that
$$
U_{R=L}\supseteq(R-2L-1)U,\eqno(1.13)$$
$$
U_{R=2L+1}\supseteq(R-L)U.\eqno(1.14)$$
Any of (1.13), (1.14) implies
$$
\dim U_{R=L}\geq \dim U-\dim U_{R=2L+1}.\eqno(1.15)$$
\par
The last inequality is equivalent to (1.12) 
provided that
$$
P=U_{R=L}\cap U_{R=2L+1}=0.\eqno(1.16)$$
\par
To prove (1.16), assume that $P\not=0.$
Then we have
$$
R=L=2L+1\ {\rm on}\ P.$$
Hence
$$
R=L=-\rm {Id}\; {\rm on}\ P.$$
Thus we come to a contradiction with Lemma 1.1.
Lemma 1.2 is proved.
\par
\vskip 3mm
{\bf Corollary of the proof.}
\quad
{\it The following formulas are true:}
$$
U_{R=L}=(R-2L-1)U,\eqno(1.17)$$
$$
U_{R=2L+1}=(R-L)U.\eqno(1.18)$$
\par
Indeed, we have proved, that we have an  equality in (1.15), which also
means that 
there are the equalities in (1.13) and (1.14).
\par
\vskip 3mm
{\bf Lemma 1.3.}\quad
{\it The subspaces $U_{R=L}$ and $U_{R=2L+1}$ are
invariant with respect to $L$ and $R$.}
\par
\vskip 3mm
{\it Proof.}\quad
Let us prove, for example,
$$
LU_{L=R}\subseteq U_{L=R}.$$
\par
Using (1.17) and (1.6),
we have
$$
LU_{L=R}=L(R-2L-1)
U=(R-2L-1)LU\subseteq U_{L=R}.$$

The other assertions are proved in the same way.
The lemma is proved.
\par
Thus we can consider the actions of $L$ and  $R$ separately on
$U_{L=R}$ and on
$U_{L=2R+1}.$
\par
\vskip 3mm
{\bf Corollary.}\quad
{\it There is no  vector  } $a \neq 0$  {\it such that} 
$$
La=-a.\eqno(1.18^{'})$$
\par
\vskip 3mm
{\it Proof.}\quad
Let $a=a_{1}+a_{2}$ where
$a_{1}\in U_{R=L}$ and
$a_{2}\in U_{R=2L+1}.$
Then it follows from Lemma 1.3 that
$$
La_{1}=-a_{1},\ La_{2}=-a_{2}.$$
\par
Considering the operator $R$ on $U_{R=L}$
and $U_{R=2L+1},$ we also have
$$
Ra_{1}=-a_{1},\ Ra_{2}=-a_{2}.$$
Now the Corollary follows from Lemma 1.1.
\par
\vskip 3mm
{\bf Definition 2.}
\quad
We shall call a subspace of all vectors $a$ satisfying
$$
La=\lambda a,\ Ra=\mu a \eqno(1.19)$$
{\it an eigen-subspace of the tripotent} $e$ 
corresponding to 
{\it the eigen-values}
$\lambda,\mu$
and denote it by
$U_{\lambda\mu}.$
\par
First of all we prove that
there are only two possibilities:
\par
1)\ $\mu=\lambda$ and 2)\ $\mu=2\lambda+1.$
\par
We will denote by $U_{\mu=\lambda}$ the sum of all
subspaces with the first possibility and by
$U_{\mu=2\lambda+1}$
with the second one.
\par
\vskip 3mm
{\bf Lemma 1.4.}\quad
\it
Let $a\in U_{\lambda \mu}.$
Then either
$a\in U_{\mu=\lambda}$ or $a\in U_{\mu=2\lambda+1}.$
\par
\vskip 3mm
{\it Proof.}\quad
\rm
Suppose
$a\notin U_{\mu=\lambda}, U_{\mu=2\lambda+1}.$
Let us present
$$
a=a_{1}+a_{2},\ a_{1}\in U_{R=L},\ a_{2}\in U_{R=2L+1},$$
where
$a_{1}\not= 0, a_{2}\not= 0.$
\par
By definition of $U_{\lambda\mu}$
$$
L(a_{1}+a_{2})=
\lambda a_{1}+\lambda a_{2},\;
R(a_{1}+a_{2})=
\mu a_{1}+\mu a_{2},$$
This implies according to Lemma 1.3
$$
La_{1}=\lambda a_{1},\;
La_{2}=\lambda a_{2},\;
Ra_{1}=\mu a_{1},\;
Ra_{2}=\mu a_{2}.$$
Hence both $a_{1}$ and $a_{2}$
are  common eigen-vectors of $L$ and
$R$.
\par
This means, that
$\mu =\lambda, \mu =2\lambda +1,$
which implies
$\lambda=-1.$
\par
But this contradicts the  Corollary of Lemma 1.3.
The lemma is proved.
\par
\vskip 3mm
{\bf Lemma 1.5.}\quad
\it
Let $a\in U_{\mu=\lambda}.$
Then one of the following three possibilities occurs:
\par
a)\ $a\in U_{00}$ or $a\in U_{{1\over 2}{1\over 2}},$
and $\;Q(a)=0$,
\par
b)\ $a\in U_{11},$ and $\;Q(a)\in U_{11},\;Q^{2}(a)=a$,
\par
c)\ $a \in U_{{3\over 2}{3\over 2}},$
and $\;Q(a)\in U_{{1\over 2}2},\;Q^{2}(a)=3a.$
\par
\vskip 3mm
{\it Proof.}
\quad
\rm
First suppose that
$Q(a)=0.$
Then (1.4) is equivalent to 
$$
(2L^{2}-L)a=2L(L-{1\over 2})a=0.\eqno(1.20)$$
Thus $\lambda = 0, {1\over 2}$ and we come to the case a).
 \par
 Let $Q(a)\not= 0.$
Acting on $a$ by both sides of equations (1.7) and (1.5)
we obtain
$$
LQ(a)=(2-\lambda )Qa,\eqno(1.21)$$
$$
(R+L)Q(a)=(1+\lambda)Qa.$$
Subtracting we have
$$
RQ(a)=(2\lambda-1)Qa.\eqno(1.22)$$
Thus $Q(a) \in U_{2-\lambda,2\lambda-1}.$
\par
According to Lemma 1.4
$Q(a)$ belongs to
$U_{\mu=\lambda}$
or to
$U_{\mu=2\lambda+1}.$
In the first case we have
$$
2\lambda -1=2-\lambda,$$
which implies $\lambda=1$.
\par
The equality $\lambda =1$ means both $a\in U_{11}$ and $Qa\in U_{11}.$
Finally acting by (1.4) on $a$ we obtain
$$
Q^{2}(a)=a.$$
Thus we come to the case b).
\par
In the second case $Q(a)\in U_{\mu=2\lambda+1}$.
Then it follows from (1.21) and (1.22) that
$$
2\lambda -1=2(2-\lambda)+1,$$
which implies
$\lambda={3\over 2}.$
Hence $a\in U_{{3\over 2}{3\over 2}},\ 
Qa\in U_{{1\over 2}2}.$
\par
Acting by (1.4) on $a$ we obtain
$$
Q^{2}(a)=3a.$$
Thus we come to the case c).
The lemma is proved.
\par
\vskip 3mm
{\bf Lemma 1.6.}\quad
\it
Let $a\in U_{\mu=2\lambda+1}.$
Then one of the following three possibilities occurs:
\par
a$'$)\ $a\in U_{01}$\ or $a\in U_{-{1\over 2}0},$
and $Q(a)=0,$
\par
b$'$)\ $a\in U_{13}, \; Q(a)\in U_{13},$
and $Q^{2}(a)=9a,$
\par
c$'$)\ $a\in U_{{1\over 2}2},\; Q(a)\in 
U_{{3\over 2}{3\over 2}},$
and $Q^{2}(a)=3a.$
\par
\vskip 3mm
{\it Proof.}\quad
\rm
First suppose that
$Q(a)=0.$
Then (1.4) is equivalent to
$$
(2L+1)^{2}a+
L(2L+1)a-
(2L+1)a=3(2L+1)La=0.\eqno(1.23)$$
This means that
$\lambda=0,-{1\over 2}.$
Thus we come to the case $a'$).
\par
Suppose now that
$Q(a)\not=0.$
Then acting by (1.7) and (1.5) on $a$, we obtain
$$
LQa=(2-\lambda)Qa,\eqno(1.24)$$
$$
(R+L)Qa=(2\lambda+2)Qa.\eqno(1.25)$$
\par
Subtracting we obtain
$$
RQ(a)=3\lambda Q(a).\eqno(1.26)$$
Thus $Q(a)\in U_{2-\lambda,3\lambda}.$
According to Lemma 1.4
$Q(a)$ belongs to
$U_{\mu=\lambda}$
or to
$U_{\mu=2\lambda+1}.$
\par
In the case
$Q(a)\in U_{\mu=2\lambda+1}$
we have
$$
3\lambda=2(2-\lambda)+1,$$
which gives
$\lambda=1.$
The equality $\lambda=1$
means both $a\in U_{13}$
and 
$Q(a)\in U_{13}.$
\par
Acting by (1.4) on $a$ we obtain
$$
Q^{2}(a)=9a.\eqno(1.27)$$
Thus we come to the case 
b$'$).
\par
Now let $Q(a)\in U_{\mu=\lambda}.$
Then comparing  (1.24) and (1.26)
we obtain
$$3\lambda=2-\lambda,$$
which gives $\lambda={1\over 2}.$
Hence
$a\in U_{{1\over 2},2},\;
Q(a)\in U_{{3\over 2}{3\over 2}}.$
Acting by (1.4) on $a$ we obtain $Q^{2}(a)=3a.$
Thus we come to the case c$'$).
The lemma is proved.
\par
\vskip 3mm
{\bf Theorem 1.1}\quad
\it
Let the three linear operators
$L,R,Q$ defined on a linear space
$U$ be a solution of the system of
equations
(1.4)-(1.8).
Then the space $U$ is a direct sum
$$
U=U_{00}\oplus U_{{1\over 2}{1\over 2}}
\oplus U_{11}\oplus
U_{{3\over 2}{3\over 2}}
\oplus
U_{-{1\over 2}0}
\oplus
U_{01}
\oplus
U_{{1\over 2}2}
\oplus
U_{13},
\eqno(1.28)$$
where
\par
1)\ for any $a\in U_{\lambda \mu}$
$$
La=\lambda a,\;\; Ra=\mu a,\eqno(1.29)$$
\par
2)\ $Q(a)=0\;\; \; \forall a\in U_{00},\;U_{{1\over 2}{1\over 2}},\;
U_{01},\;U_{-{1\over 2}0},$
\par
3)\ the subspaces $U_{11}$ and $U_{13}$
are direct sums
$$
U_{11}=U_{11}^{+}
\oplus
U_{11}^{-},\ \;   U_{13}=U_{13}^{+}
\oplus
U_{13}^{-}\eqno(1.30)$$
and
$$
Q(a)=\pm a\ \; {\rm if}\ a\in U_{11}^{\pm},\eqno(1.31)$$
$$
Q(a)=\pm 3a\ \;{\rm if}\ a\in U_{13}^{\pm},\eqno(1.32)$$
\par
4)\ there is a one-to-one correspondence
$\tau=\tau^{-1}$ between the subspaces
$U_{{3\over 2}{3\over 2}}$ and
$U_{{1\over 2}2}$
such that
$$
Q(a)=\sqrt{3}\ \tau(a),\; \forall a\in U_{{3\over 2}{3\over 2}},\
U_{{1\over 2}2}.\eqno(1.33)$$
\par
Conversely, let a space $U$ be presented as a direct sum (1.28) and
three operators $L,R,Q$ on $U$ be defined by
properties 1)-4).
Then the system of equations
(1.4)-(1.8) is fulfilled.
\par
\vskip 3mm
{\it Proof.}\quad
\rm
To prove (1.28) 
and also properties 1), 2) it is enough to prove 
$$U_{R=L}=U_{\mu=\lambda}\ 
{\rm and}\ U_{R=2L+1}=U_{\mu=2\lambda+1}.\eqno(1.34)$$
Then (1.28) and properties 1), 2)
will follow from Lemmas 1.5 and
1.6.
\par
The equalities (1.34)
mean
that the operator $L$
has no Jordan blocks of second degree, i.e.,
there are no vectors
$a_{1},a_{2}$
such that
$$
La_{1}=\lambda a_{1},\ 
La_{2}=\lambda a_{2}+a_{1}.$$
\par
According to Lemmas 1.2 and 1.3 we can consider separately the two cases
$a_{1},a_{2}\in U_{R=L}$
and $a_{1},a_{2}\in U_{R=2L+1}$.
\par
In the case where
$a_{1},a_{2}\in U_{R=L}$
we have
$$
Ra_{1}=\lambda a_{1},\ Ra_{2}=\lambda a_{2}+a_{1}.\eqno(1.35)$$
\par
It follows from (1.7) that
$$
LQ(a_{1})=(2-\lambda)Q(a_{1}),$$
$$
LQ(a_{2})=(2-\lambda)Q(a_{2})-Q(a_{1}).\eqno(1.36)$$
\par
Using (1.36) we obtain from (1.5):
$$
RQ(a_{1})=(2\lambda-1)Q(a_{1}),$$
$$
RQ(a_{2})=
(2\lambda-1)Q(a_{2})+2Q(a_{1}).\eqno(1.37)$$
Equalities (1.36), (1.37) imply
$$
(R-2L-1)Qa_{1}=(4\lambda-6)Qa_{1},$$
$$
(R-2L-1)Qa_{2}=(4\lambda-6)Qa_{2}+4a_{1},$$
$$
(R-L)Qa_{2}=3(\lambda-1)Qa_{2}+3Qa_{1}.$$
\par
Now acting by (1.8) on $Q(a_{2})$ and
using the equalities above we obtain
$$
12(\lambda-1)(\lambda-{3\over 2})Q(a_{2})+
12(2\lambda-{5\over 2})Qa_{1}=0.
\eqno(1.38)$$
\par
According to Lemma 1.5,
$\; \lambda=0,{1\over 2},1,{3\over 2},$
and $Q(a_{1})\not= 0,$
if
$\lambda=1,{3\over 2}.$
Thus in the case where $\lambda=1,{3\over 2}$  the
first term is equal to zero, but not the
second one.
Hence these  cases 
are impossible.
\par
In the remaining cases 
  $\lambda=0,{1\over 2}$ we have $Q(a_{1})=0$ but
it follows from (1.38)
that not only $Q(a_{1})$
but also $Q(a_{2})$ is equal to zero.
\par
Let us apply in these cases the relation (1.4) to $a_{2}.$
We obtain
$$
(2\lambda-1)\lambda a_{2}+
(4\lambda-1)a_{1}=0.$$
Hence in both cases
($\lambda=0,{1\over 2}$)
we come to contradiction with $a_{1}\neq 0$.
\par
In the same way we consider the case where
$a_{1},a_{2}\in U_{R=2L+1}$.
Then the relation (1.36) is valid, but instead of (1.35) we have
$$
Ra_{1}=(2\lambda+1)a_{1},\ 
Ra_{2}=(2\lambda+1)a_{2}+2a_{1},\eqno(1.35')$$
and instead of (1.37) we have
$$
RQa_{1}=3\lambda Qa_{1},$$
$$
RQa_{2}=3\lambda Qa_{2}+3Qa_{1}.\eqno(1.37')$$
\par
Equalities (1.36), (1.37')  imply
$$
(R-2L-1)Qa_{1}=(5\lambda-5)Qa_{1},$$
$$
(R-2L-1)Qa_{2}=(5\lambda-5)Qa_{2}+5Qa_{1},$$
$$
(R-L)Qa_{2}=(4\lambda+2)Qa_{2}+4Qa_{1}.$$
Again acting by (1.8) on $Q(a_{2})$ and using the
equalities above we obtain
$$
20(\lambda-1)(\lambda-{1\over 2})Qa_{2}
+20(2\lambda-{1\over 2})Qa_{1}=0.
\eqno(1.38')$$
\par
According to Lemma 1.6,
$\; \lambda=-{1\over 2},0,{1\over 2},1$
and
$Q(a_{1})\not=0$
if $\lambda={1\over 2},1.$
Thus in the case  where $\lambda={1\over 2},1$
 the first term is equal to zero but not the second one.
Hence these  cases  are impossible.
\par
In the remaining cases
($\lambda=0,-{1\over 2}$) we have  $Q(a_{1})=0$. But
it follows from (1.38')
that not only $Q(a_{1}),$ but also $Q(a_{2})$
is equal to zero.
Let us act in these cases by relation (1.4) on $a_{2}$.
We obtain
$$
3\lambda(2\lambda+1)a_{2}+
3(5\lambda+1)a_{1}=0.$$
\par
Hence
in both cases $\lambda=0,-{1\over 2}$
we come to a contradiction with $a_{1}\neq 0$ .
Thus the equality (1.28)
and properties 1), 2) are proved.
\par
The property 3) for $U_{11}$ and
$U_{13}$
is an elementary consequence of the properties
$$
Q(a)\in U_{11},\ 
Q^{2}(a)=a\ {\rm for}\ a\in U_{11},$$
$$
Q(a)\in U_{13},\ Q^{2}(a)=9a\ {\rm for}\ a\in U_{13},$$
which were proved in Lemmas 1.5, 1.6.
\par
To prove 4) we just define
$$
\tau a={1\over \sqrt{3}}Qa,\ \forall a\in
U_{{3\over 2}{3\over 2}},\;
U_{{1\over 2}2}.\eqno(1.39)$$
Then the correctness of this definition and the property
$\tau^{-1}=\tau$ follow from
$Q^{2}(a)=3a.$
\par
The converse assertion can be checked directly for all 
ten subspaces in
(1.28), (1.30).
The theorem is proved.
\par
\vskip 3mm
{\bf Definition 3. ([13])}
\quad
A generalized Jordan triple system of second order
 $U$ is said to be 
{\it weakly commutative} if 
$$
(yx(yyy))=((yyy)xy),\ \forall x,y\in U.
\eqno(1.40)$$
\par
We will consider an example of a weakly commutative triple
system in \S 2.
\par
The identity (1.40)
implies an additional condition on the operators
$L,Q,R.$
To find it consider the polarization of (1.40):
$$
(ux(yyy))+(yx(uyy))+
(yx(yuy))+(yx(yyu))=$$
$$
=((yyy)xu)+((uyy)xy)+
((yuy)xy)+((yyu)xy).\eqno(1.41)$$
\par
If we substitute $y=x=e$,
we obtain the equation
$$
(R-L)(R+Q+L-1)=0.\eqno(1.42)$$
\par
Subtracting this equation from (1.8$'$), we get
$$
(R-L)Q=-3(R-L)L.\eqno(1.43)$$
\par
\vskip 3mm
{\bf Theorem 1.2.}\quad
{\it 
Let three linear operators $L,\ R,\ Q$
defined on a linear space $U$ be
a solution of the system of equations
(1.4)-(1.8),\ (1.42).
Then the space $U$ 
is a direct sum of six subspaces
$U_{\lambda\mu}$:
$$
U=U_{00}\oplus U_{01}
\oplus U_{{1\over 2}{1\over 2}}
\oplus
U_{11}^{+}\oplus 
U_{11}^{-}\oplus 
U_{13}^{-},\eqno(1.44)$$
such that
for every $a\in U_{\lambda\mu}$
$$
La=\lambda a,Ra=\mu a,\eqno(1.45)$$
and
$$
Qa=0,\ {\rm if}\ \lambda\not=1,\eqno(1.46)$$
$$
Qa=-3a\ {\rm if }\ a\in U_{13}^{-},\eqno(1.47)$$
$$
Qa=\pm a,\forall  a\in U_{11}^{\pm} .\eqno(1.48)$$

Conversely,
let a space $U$ be presented as a direct sum (1.44)
and three operators $L,R,Q$ on $U$ be defined by
formulas (1.45)-(1.48).
Then the system of equations (1.4)-(1.8),\ (1.42)
is fulfilled.
\par
{\it Proof.}\quad
\rm
It follows from Theorem 1.1 that the space $U$
is in general a direct sum of
ten subspaces.
We have to prove that the additional property of weak
commutativity implies
that four of them
$
U_{{3\over 2}{3\over 2}},\ U_{-{1\over 2}0},\
\ U_{{1\over 2}2},\ U_{13}^{+}\ {\rm are\ zero.}$
\par
We rewrite (1.43) as
$$
(R-L)Qa=-3(R-L)La.\eqno(1.49)$$
\par
First suppose that $a\in U_{{3\over 2}{3\over 2}}.$
Then the right-hand side equals  zero.
But the left-hand side is not equal to zero,
because $\tau=$$1\over \sqrt{3}$$ Q$ defines a one to-one-correspondence between
 $ U_{{3\over 2}{3\over 2}}$ and $ U_{{1\over 2}2}$ which implies $ Q(a)\neq 0$ 
and $ Q(a)\in U_{{1\over 2}2}. $ 
 \par
Now suppose that
$a\in U_{-{1\over 2}0},\; U_{{1\over 2}2}.$
Then the left-hand side equals zero but the right 
hand-side is not equal to 
 zero as  is easy to check.
\par
Finally let $a\in U_{13}^{+}.$
Then the left-hand side equals $6a$ but the right-hand side equals $-6a.$
\par
Thus in all 4 cases we come to a contradiction with $a\not=0.$
Hence all 4  considered subspaces are zero.

The converse assertion can be checked directly for all
six subspaces in (1.44).
The theorem is proved.
\par
The next theorem is an evident corollary of Theorems 1.1 and 1.2.
\par
\vskip 3mm
{\bf Theorem 1.3.}\quad
\it Let $e$ be a tripotent of a generalized Jordan 
triple system of second order
and three operators $L,R,Q$ be defined by the formulas (1.3).
Then the space $U$ of the triple system is a direct sum of ten
subspaces
$$
U=U_{00}\oplus U_{{1\over 2}{1\over 2}}\oplus
U_{11}^{+}\oplus U_{11}^{-}\oplus
U_{{3\over 2}{3\over 2}}\oplus
U_{-{1\over 2}0}\oplus
U_{01}\oplus U_{{1\over 2}2}\oplus
U_{13}^{+}\oplus U_{13}^{-},\eqno(1.50)$$
and the operators $L,R,Q$ have the following properties:
$$
1)\ {\it for\ any}\ a\in U_{\lambda\mu},\ La=\lambda a,\;
Ra=\mu a,\eqno(1.51)$$
$$
2)\ Q(a)=0\ \forall a\in U_{00},\; U_{{1\over 2}{1\over 2}},\;
U_{01},\;
U_{{1\over 2}0},\eqno(1.52)$$
$$
3)\ Q(a)=\pm a,\ if\ a\in U_{11}^{\pm},\;\;
Q(a)=\pm 3a,\ if\ a\in U_{13}^{\pm},\eqno(1.53)$$ 
$$
4)\ {\it the\ map}\ 
\tau a={1\over \sqrt{3}}Q(a)\ \forall
a\in U_{{3\over 2}{3\over 2}},\;U_{{1\over 2}2}\eqno(1.54)$$
is a one-to-one correspondence between the subspaces
$U_{{3\over 2}{3\over 2}}$
and $U_{{1\over 2}2}.$
\par
If in addition the triple system is weakly commutative, then
the space $U$ is a direct sum of only  six subspaces
$$
U=U_{00}\oplus U_{{1\over 2}{1\over 2}}\oplus
U_{11}^{+}\oplus U_{11}^{-}\oplus
U_{01}\oplus 
U_{13}^{-},\eqno(1.55)$$
and the operators $L,Q,R$
have the  properties (1.51),\ (1.52),\ (1.53).
\par
\vskip 3mm
{\bf Definition 4.}\quad
\rm
We shall call the direct sum (1.50)
{\it the Peirce decomposition defined by a tripotent $e$}.
\par
\vskip 3mm
\begin{center}
{\bf \S\ 2\ Examples of the Peirce decomposition}
\end{center}
\par
In this section, we will give several examples 
of the Peirce decomposition for
generalized Jordan triple systems of second 
order.
\par
\vskip 3mm
{\bf 1. Triple system $A_{kn}-A_{nk}.$}
\par
The space of the triple is the set of pairs
$\left(\matrix{A_{1}\cr A_{2}\cr}\right),$
where $A_{1}$
is a $(k,n)$-matrix and $A_{2}$ is a $(n,k)$-matrix.
The triple product is given by the formula:
$$
\left(\left(\matrix{A_{1}\cr A_{2}\cr}\right)
\left(\matrix{B_{1}\cr B_{2}\cr}\right)
\left(\matrix{C_{1}\cr C_{2}\cr}\right)\right)=
\left(\matrix{A_{1}B_{1}^{T}C_{1}+C_{1}B_{1}^{T}A_{1}-C_{1}A_{2}B_{2}^{T}
\vspace{1mm}\cr 
A_{2}B_{2}^{T}C_{2}+C_{2}B_{2}^{T}A_{2}-B_{1}^{T}A_{1}C_{2}\cr}\right).\eqno(2.1)$$
Let us consider the case $k<n$ and present the
matrices $A_{1},A_{2}$
in the form
$A_{1}=(A_{11},A_{12})$
where $A_{11}$ is a $(k,k)$-matrix,
$A_{12}$ is a
$(k,n-k)$-matrix and
$A_{2}=(A_{21},A_{22})$,
where $A_{21}$ is a
$(k,k)$-matrix,
$A_{22}$ is a
$(n-k,k)$-matrix.
Thus the elements of the space will be
presented as quadruples
$$\left(
\begin{array}{c|c}
A_{11}&A_{12}\\
\hline
A_{21}&A_{22}\\
\end{array}\right).$$
\par
It is easy to cheek that the element
$$e=\left(
\begin{tabular}{c|c}
$E_{k}$&$0$\\
\hline
$E_{k}$&0\\
\end{tabular}\right),
$$
where $E_{k}$ is the identity matrix of order
$k$ is a tripotent.
The operators
$L,R,Q$
in this case act as follows
$$
L\left(
\begin{tabular}{c|c}
$A_{11}$&$A_{12}$\\
\hline
$A_{21}$&$A_{22}$\\
\end{tabular}\right)
=\left(
\begin{tabular}{c|c}
$A_{11}$&$A_{12}$\\
\hline
$A_{21}$&$A_{22}$\\
\end{tabular}\right),
\eqno(2.2)
$$
$$
R\left(
\begin{tabular}{c|c}
$A_{11}$&$A_{12}$\\
\hline
$A_{21}$&$A_{22}$\\
\end{tabular}\right)
=
\left(
\begin{tabular}{c|c}
$2A_{11}-A_{21}$&$A_{12}$\\
\hline
$2A_{21}-A_{11}$&$A_{22}$\\
\end{tabular}\right),
\eqno(2.3)$$
$$
Q\left(
\begin{tabular}{c|c}
$A_{11}$&$A_{12}$\\
\hline
$A_{21}$&$A_{22}$\\
\end{tabular}\right)
=
\left(
\begin{tabular}{c|c}
$2A_{11}^{T}-A_{21}^{T}$&$-A_{22}^{T}$\\
\hline
$2A_{21}^{T}-A_{11}^{T}$&$-A_{21}$\\
\end{tabular}\right).\eqno(2.4)$$
Thus
the Peirce decomposition has the form
$$
U=U_{11}^{+}\oplus
U_{11}^{-}\oplus U_{13}^{+}\oplus U_{13}^{-},\eqno(2.5)$$
where
$$
U_{11}^{+}=\left\{
\left(
\begin{tabular}{c|c}
$S_1$&$A$\\
\hline
$S_1$&$-A^{T}$\\
\end{tabular}\right)\right\},\ 
U_{11}^{-}=\left\{
\left(
\begin{tabular}{c|c}
$K_1$&$B$\\
\hline
$K_1$&$B^{T}$\\
\end{tabular}\right)\right\},$$
$$
U_{13}^{+}=\left\{
\left(
\begin{tabular}{c|c}
$S_2$&$0$\\
\hline
$-S_2$&$0$\\
\end{tabular}\right)\right\},\ 
U_{13}^{-}=\left\{
\left(
\begin{tabular}{c|c}
$K_2$&$0$\\
\hline
$-K_2$&$0$\\
\end{tabular}\right)\right\},$$
where $A,B$ are arbitrary $(k,n)$ matrices, $S_1,S_2$ (respectively
 $K_1,K_2$ ) are arbitrary symmetric (respectively
 skew-symmetric ) matrices of order $k$.

In the next example we consider a tripotent,
such that all ten components in the
Peirce decomposition are nonzero.
\par
\vskip 3mm
{\bf 2. Triple $A_{nn}-A_{nn}.$}
\par
This triple system is a special case of the previous triple
with $k=n$.
Thus the space consists of pairs of matrices
$\left(\matrix{A_{1}\cr A_{2}\cr}\right)$
of order $n$ and the product is given by (2.1).
We will consider the  case $n=3l$ and present matrices
$A_{1},A_{2}$ as block matrices of order
$3$ with matrices of order $l$ as elements.
\par
It is easy to check that the
element
$e=\left(\matrix{e_{1}\cr e_{2}\cr}\right),$
where
$$
e_{1}=
\left(
\begin{tabular}{c|c|c}
${1\over \sqrt{2} }E_{l}$&0&0\\
\hline
0&$E_{l}$&0\\
\hline
0&0&0\\
\end{tabular}\right),\ 
e_{2}=
\left(
\begin{tabular}{c|c|c}
0&0&0\\
\hline
0&$E_{l}$&0\\
\hline
0&0&${1\over \sqrt{2}}E_{l}$\\
\end{tabular}\right)
\eqno(2.6)$$
and $E_{l}$ is the identity matrix
of order $l$, is a tripotent.
\par
We denote the elements of the  pair
$\left(\matrix{X\cr Y\cr}\right)$
by
$$
X=
\left(
\begin{tabular}{c|c|c}
$X_{11}$&$X_{12}$&$X_{13}$\\
\hline
$X_{21}$&$X_{22}$&$X_{23}$\\
\hline
$X_{31}$&$X_{32}$&$X_{33}$\\
\end{tabular}\right),\ 
Y=\left(
\begin{tabular}{c|c|c}
$Y_{11}$&$Y_{12}$&$Y_{13}$\\
\hline
$Y_{21}$&$Y_{22}$&$Y_{23}$\\
\hline
$Y_{31}$&$Y_{32}$&$Y_{33}$\\
\end{tabular}\right).$$
Then the action of operators
$L,R,Q$
is given by the  formula
$$
L\left(\matrix{X\cr Y\cr}\right)=
\left(\matrix{
\left(
\begin{tabular}{c|c|c}
$X_{11}$&${1\over 2}X_{12}$&0\\
\hline
${3\over 2}X_{21}$&$X_{22}$&${1\over 2}X_{23}$\\
\hline
${1\over 2}X_{31}$&0&$-{1\over 2}X_{33}$\\
\end{tabular}\right)\vspace {2mm} \cr\\  
\left(
\begin{tabular}{c|c|c}
$-{1\over 2}Y_{11}$&${1\over }Y_{12}$&0\\
\hline
$0$&$Y_{22}$&${1\over 2}Y_{23}$\\
\hline
${1\over 2}Y_{31}$&${3\over 2}Y_{32}$&$Y_{33}$\\
\end{tabular}\right)\cr}\right),
\eqno(2.7)$$

$$
R\left(\matrix{X\cr Y\cr}\right)=
\left(\matrix{\left(
\begin{tabular}{c|c|c}
$X_{11}$&${3\over 2}X_{12}-{1\over \sqrt{2}}Y_{12}$&
${1\over 2}X_{13}-{1\over 2}Y_{13}$\\
\hline
${3\over 2}X_{21}$&$2X_{22}-Y_{22}$&
$X_{23}-{1\over \sqrt{2}}Y_{23}$\\
\hline
${1\over 2}X_{31}$&$X_{32}$&$0$\\
\end{tabular}\right)\vspace {2mm} \cr\\ 
\left(
\begin{tabular}{c|c|c}
$0$&$Y_{12}-{1\over \sqrt{2}}X_{12}$&
${1\over 2}Y_{13}-{1\over 2}X_{13}$\\
\hline
$Y_{21}$&$2Y_{22}-X_{22}$&${3\over 2}Y_{23}-{1\over {\sqrt 2}}X_{23}$\\
\hline
${1\over 2}Y_{31}$&${3\over 2}Y_{32}$&$Y_{33}$\\
\end{tabular}\right)\cr}\right),
\eqno(2.8)$$
$$
Q\left(\matrix{X\cr Y\cr}\right)=
\left(\matrix{
\left(
\begin{tabular}{c|c|c}
$X_{11}^{T}$&$\sqrt{2}X_{21}^{T}$&0\\
\hline
$\sqrt{2}X_{12}^{T}-Y_{12}^{T}$&$2X_{22}^{T}-Y_{22}^{T}$&$-Y_{32}^{T}$\\
\hline
$0$&0&$0$\\
\end{tabular}\right)\vspace {1mm}\cr \vspace {1mm}
\left(
\begin{tabular}{c|c|c}
$0$&$-X_{21}^{T}$&0\\
\hline
$0$&$2Y_{22}^{T}-X_{22}^{T}$&$\sqrt{2}Y_{32}^{T}$\\
\hline
$0$&$\sqrt{2}Y_{23}^{T}-X_{23}^{T}$&$Y_{33}^{T}$\\
\end{tabular}\right)\cr}\right).
\eqno(2.9)$$
The Peirce decomposition for this tripotent has the form
(1.50),
where all ten components are nonzero.
They have the following form
(we define only elements which are  not necessarily zero)
$$
U_{00}=\{X_{13}=Y_{13}=A_1\},$$
$$ 
U_{{1\over 2}{1\over 2}}=\{Y_{12}=\sqrt{2}X_{12}=\sqrt{2}B_1,\ 
X_{23}=\sqrt{2}Y_{23}=\sqrt{2}C_1,\ X_{31}=D,\;Y_{31}=E\},$$
$$
U_{11}^{+}=\{X_{11}=S_{1},\;X_{22}=Y_{22}=S_{2},\;
Y_{33}=S_{3}\},$$
$$
U_{11}^{-}=\{X_{11}=K_{1},\;X_{22}=Y_{22}=K_{2},\;Y_{33}=K_{3}\},$$
$$
U_{{3\over 2}{3\over 2}}=\{X_{21}=F,\;Y_{32}=G\},$$
$$
U_{-{1\over 2}0}=\{X_{33}=H,\;Y_{11}=I\},$$
$$
U_{01}=\{X_{13}=-Y_{13}=A_2,\;X_{32}=J,\;Y_{23}=L\},$$
$$
U_{{1\over 2}2}=\{X_{12}=-\sqrt{2}Y_{12}=\sqrt{2}B_2,\;
Y_{32}=-\sqrt{2}X_{32}=\sqrt{2}C_2\},$$
$$
U_{13}^{+}=\{X_{22}=S_4,\;Y_{22}=-S_4\},\;\;
U_{13}^{-}=\{X_{22}=K_4,\;Y_{22}=K_4\},$$
where
$A_1,A_2,B_1,B_2,C_1,C_2,D,E,F,G,H,I,J,L$
are arbitrary matrices of order $l;\;\;$
 $S_{1},S_{2},S_{3},S_{4}$
are arbitrary symmetric matrices of order $l;\;\;$
 $K_{1},K_{2},K_{3},K_{4}$
are arbitrary skew-symmetric matrices of order $l$.
\par
In the previous two examples the triple 
systems were not weakly commutative.
In the next one the triple system will have this property. 
\par 
\vskip 3mm
{\bf 3.\ Triple system $D_{nk}$}
\par
\vskip 3mm
The triple system is the space of all $n\times k$
matrices,
with the
triple product defined by
$$
(XYZ):=
XY^{T}Z+ZY^{T}X-YX^{T}Z,\eqno(2.10)$$
where $X^{T}$
denotes the transpose of $X$.
\par
The triple system $D_{nk}$ is weakly commutative.
Indeed,
$$
(XY(XXX))-((XXX)YX)=-YX^{T}XX^{T}X+
Y(XX^{T}X)^{T}X=0.$$
Thus according to the Theorem 1.2 we can not have more than six components 
in the Peirce decomposition.
\par
We will consider the case where $n\geq k$.
\par
It is easy to check that the following element $e$ is a tripotent:
$$
e=\left(
\begin{tabular}{c|c}
E&0\\
\hline
0&0\\
\end{tabular}\right),\eqno(2.11)
$$
where $E$ is the identity matrix of order
$l\le k.$
\par
Let
$$
X=\left(
\begin{tabular}{c|c}
A&B\\
\hline
C&D\\
\end{tabular}\right)
$$
be an arbitrary element,
where
$A$ is an $l\times l$ matrix, $B$ is an $l\times(k-l)$ matrix,
$C$ is an $(n-l)\times l$ matrix, $D$ is an $(n-l)\times (k-l)$ matrix.
\par
Then we have
$$
L(X)=\left(
\begin{tabular}{c|c}
A&0\\
\hline
C&0\\
\end{tabular}\right),$$
$$
R(X)=\left(
\begin{tabular}{c|c}
$2A-A^{T}$&B\\
\hline
C&0\\
\end{tabular}\right),
$$
$$
Q(X)=\left(
\begin{tabular}{c|c}
$2A^{T}-A$&0\\
\hline
-C&0\\
\end{tabular}\right).$$
\par
Thus, in this case
we obtain
$$
U=U_{00}\oplus U_{11}^{+}\oplus
U_{11}^{-}\oplus U_{01}\oplus
U_{13}^{-},\eqno(2.12)$$
where the five subspaces of matrices
are as follows:
$$U_{00}=\left(
\begin{tabular}{c|c}
0&0\\
\hline
0&D\\
\end{tabular}
\right),\;
U_{11}^{+}=\left(
\begin{tabular}{c|c}
$A_{s}$&$0$\\
\hline
$0$&$0$\\
\end{tabular}
\right),\  
U_{11}^{-}=\left(
\begin{tabular}{c|c}
0&0\\
\hline
C&0\\
\end{tabular}
\right),
$$
$$
U_{01}=\left(
\begin{tabular}{c|c}
0&B\\
\hline
0&0\\
\end{tabular}
\right),\ 
U_{13}^{-}=\left(
\begin{tabular}{c|c}
$A_{k}$&0\\
\hline
0&0\\
\end{tabular}
\right).\   
$$
Here by $A_{s}$ and $A_{k}$ we denote the symmetric and
skew-symmetric parts of the 
matrix $A$:
$$
A_{s}={1\over 2}(A+A^{T}),\ A_{k}={1\over 2}(A-A^{T}).$$
We note that in the special case where
$k=l$,
one has $Lx=x$
for all $x$ and the Peirce decomposition has the form
$$
U=U_{11}^{+}\oplus U_{11}^{-}\oplus U_{13}^{-}.\eqno(2.13)$$
\par
\vskip 3mm
{\bf 4.\ 
Generalized Jordan triple systems of
second order
 defined by 
structurable algebras}
\par
\vskip 3mm
By definition ([1]) the structurable algebra is
an algebra with a multiplication
$x\circ y$ and involutive antiautomorphism  $x \rightarrow \bar x$ such that 
the triple product
$$
(xyz)=(x\circ {\bar y})\circ
z+(z\circ {\bar y})\circ x
-(z\circ{\bar x})\circ y
\eqno(2.14)
$$
defines a 
generalized Jordan triple system of
second order.
\par
\vskip 3mm
Let $e$ be a hermitian idempotent  of a structurable algebra,
i.e.
$e\circ e=e,\;
{\bar e}=e.$
Clearly $e$ is a tripotent.
\par
Such tripotents with one more additional condition were studied by
B. Allison ([2]). We give as an example the Peirce decomposition in  the
simplest case of such  tripotent, where $e$ is the unit of the
structurable algebra.
We have 
$$
L(x)=x,R(x)=2x-{\bar x},Q(x)=2{\bar x}-x.$$
The Peirce decomposition consists only of two components:
$$U=U_{11}^{+}\oplus U_{13}^{-},\eqno(2.15)$$
where
$$
U_{11}^{+}=\{x|{\bar x}=x\},\ 
U_{13}=\{x|{\bar x}=-x\}.$$
This is the  decomposition of the space $U$ into subspaces of
symmetric and skew-symmetric elements:
$x={x+{\bar x}\over 2}+
{x-{\bar x}\over 2},$ $x \in U.$
\par

\par
\vskip 3mm
\begin{center}
\bf \S 3 Relations between subspaces in the
Peirce decomposition
\end{center}
\par
\vskip 3mm
In \S 1 we considered a system of all linear equations,
obtained from the main identities
(1.1), (1.2) by all possible four times
 substitutions of a given tripotent
$e$.
\par
Theorems 1.1 and 1.2 give the  general solution of
this system.
\par
In this section we go further and consider the system of all
bilinear equations obtained from the main identities
(1.1), (1.2)
by all possible three times
substitutions of the tripotent $e$.
This system defines bilinear relations
between subspaces in the Peirce decomposition.
We solve this system in the special but important case,
where the transformation $L$ is the identity operator.
\par
We write all possible relations between two elements.
For this we substitute in
(1.1) and (1.2)
three elements equal to $e$ in all possible combinations.
 From the equation (1.1) we obtain
\par
\vskip 2mm
\begin{center}
\begin{tabular}{l|lr}
$y,z$&$L(yez)=(L(y)ez)-(yez)+(yeL(z)),$&(3.1)\\
\hline
$y,v$&$L(yve)=(L(y)ve)-(yL(v)e)+(yve),$&(3.2)\\
\hline
$v,z$&$L(evz)=(evz)-(eL(v)z)+(evL(z)),$&(3.3)\\
\hline
$x,u$&$(xue)=R(xue)-Q(uxe)+L(xue),$&(3.4)\\
\hline
$x,v$&$(xeQ(v))=(R(x)ve)-Q(exv)+(evR(x)),$&(3.5)\\
\hline
$x,y$&$(xeR(y))=R(xey)-(yQ(x)e)+(yeR(x)),$&(3.6)\\
\hline
$x,z$&$(xeL(z))=(R(x)ez)-(eQ(x)z)+L(xez),$&(3.7)\\
\hline
$u,y$&$(euR(y))=R(euy)-(yR(u)e)+(yeQ(u)),$&(3.8)\\
\hline
$u,v$&$(euQ(v))=(Q(u)ve)-Q(uev)+(evQ(u)),$&(3.9)\\
\hline
$u,z$&$(euL(z))=(Q(u)ez)-(eR(u)z)+L(euz).$&(3.10)\\
\end{tabular}
\end{center}
\par
 From the equation (1.2) we obtain
\par
\vskip 2mm
\begin{center}
\begin{tabular}{l|lr}
$a,b$&$(R-2L-1)((aeb)-(bea))=0,
$&(3.11)\\
\hline
$a,u$&$((R-2L-1)a,u,e)-(e,u,(R-2L-1)a)=0,$&(3.12)\\
\hline
$a,v$&$(R-L)((ave)-(eva))+(eR(v)a)-(aR(v)e)=
(ev(R-L)a),$&(3.13)\\
\hline
$a,c$&$((R-L)aec)-2(ce(R-L)a)=(aQ(c)e)-(eQ(c)a).$&(3.14)\\
\end{tabular}
\end{center}
\par
We note,
that to obtain the formula
(3.12) we used formulas (3.2)
and (3.3) after the substitution.

In general on can consider the following problem .
Let a decomposition of a linear space $U$  be given
by formulas (1.28), (1.30)
and three operators 
$L,R,Q$
be defined by (1.29), (1.31), (1.32), (1.33).
Then one can consider the system of identities (3.1)-(3.14) as a system
of equations with tree unknown bilinear operators:
$(exy),(xey),(xye).$
\par
We are going to solve this system in the case where
$e$ is a left unit.
\par

\par
\vskip 3mm
{\bf Definition 5.}
\quad
A tripotent $e$ is called a {\it left unit}
if
$$
L(x)=(eex)=x\ \forall x.\eqno(3.15)$$
We note that for the left unit $e$ the Peirce decomposition 
(1.50) has the form
$$
U=U_{11}^{+}\oplus U_{11}^{-}\oplus U_{13}^{+}\oplus U_{13}^{-}.\eqno(3.16)$$
\par
According to Theorem 1.1, the 
two  transformations $R$ and $Q$  commute and they are nondegenerate
in this case.
Namely,
$$
Rx=
\left\{\matrix{x\quad {\rm if}\ x\in U_{11},\cr
3x\quad {\rm if}\ x\in U_{13},\cr}\right.\eqno(3.17)
$$
and
$$
Qx=\left\{\matrix {\pm x\quad {\rm if}\ x\in U_{11}^{\pm},\cr
\pm 3x\quad {\rm if}\ x\in U_{13}^{\pm}.\cr}\right.\eqno(3.18)$$
\par
The transformation
$x\rightarrow {\bar x},$
where 
$$
{\bar x}=RQ^{-1}x=QR^{-1}x$$
is involutive and can be defined as
$$
{\bar x}=
\left\{
\begin{array}{cl}
x & {\rm if}\ x\in U_{11}^{+},U_{13}^{+},\cr
-x & {\rm if}\ x\in U_{11}^{-},U_{13}^{-}.
\end{array}
\right.\eqno(3.19)$$
We shall call this transformation {\it the conjugation}.
\par
We will also consider another involutive transformation 
$\sim$ of the space $U$ defined by
$$
{\tilde x}=
\left\{
\begin{array}{cl}
x & {\rm if}\ x\in U_{11}^{+},\,U_{13}^{-},\cr
-x & {\rm if}\ x\in U_{11}^{-},\,U_{13}^{+}.
\end{array}
\right.\eqno(3.20)$$
\par
Moreover, we shall use the scalar function $p(x)$
on the space $U$ defined by
$$
p(x)=
\left\{
\begin{array}{cl}
0 & {\rm if}\ x\in U_{11},\cr
1 & {\rm if}\ x\in U_{13}.
\end{array}
\right.\eqno(3.21)
$$
There is a connection between the operations defined in
(3.19), (3.20) and (3.21).
It is easy to check that
$$
{\bar x}=(-1)^{p(x)}\ {\tilde x}.\eqno(3.22)$$

\par
Now we start to solve the system of equations (3.1)-(3.14) in the case where
$e$ is a left unit.
\par
First of all we note that in this case the bilinear operators
$(exy)$ and $(xye)$ could be expressed in terms of 
$(xey)$.
Even more,
the whole triple system
$(xyz)$ could be expressed in terms of  the operation
$(xey)$,
which we denote by
$$
(xey)=x\circ y.\eqno(3.23)$$
\par
Indeed,
put $b=d=e$
in the identity (1.1).
Then one obtains
$$
a\circ(c\circ f)=
(a\circ c)\circ f-(cQ(a)f)+c\circ(a\circ f).\eqno(3.24)$$
The transformation $Q$ is nondegenerate.
Thus  changing notation one can rewrite (3.24) as follows
$$
(xyz)=(Q^{-1}(y)\circ x)\circ z+
x\circ(Q^{-1}(y)\circ z)-
Q^{-1}(y)\circ(x\circ z).\eqno(3.25)$$
Moreover,
the formula (3.25) implies
$$
(exy)={\bar x}\circ y,\eqno(3.26)$$
$$
(xye)=R(Q^{-1}(y)\circ x)+x\circ {\bar y}-
Q^{-1}(y)\circ(Rx).\eqno(3.27)$$
Thus  instead of three unknowns 
$(xey),(exy),(xye)$ one can consider
only one unknown 
$(xey)=x\circ y.$
\par
\vskip 3mm
{\bf Lemma 3.1.}
\quad
\it
Let $L$ be the identity operator.
Then  formula (3.14) under conditions (3.26), (3.27) is
equivalent to the formula
$$
R(u\circ v)=
\left\{
\begin{array}{ll} 
2u\circ v-v\circ u & {\rm if}\ u,v\in U_{11},\cr
v\circ u           & {\rm if}\ u\in U_{11},v\in U_{13},\cr
4u\circ v-3v\circ u & {\rm if}\ u\in U_{13},v\in U_{11},\cr
2u\circ v-v\circ u &  {\rm if}\ u,v\in U_{13}.
\end{array}
\right.\eqno(3.28)
$$
\par
\vskip 3mm
Proof.\quad
\rm
We rewrite (3.14) using (3.26) and (3.27):
$$
R(c\circ a)=Rc\circ a-a\circ Rc+(R-1)a\circ c-c\circ (R-2)a.\eqno(3.29)$$
\par
Substituting $c=u,a=v$
it is easy to check that for all possible four cases where
$u,v\in U_{11},U_{13}$
 formula
(3.29)
becomes (3.28).
\par
It is evident that one can do these considerations backwards.
The lemma is proved.
\par
\vskip 3mm
{\bf Lemma 3.2.}\quad
\it
Let $L$ be the identity operator.
Then  formula (3.9) under conditions
(3.26), (3.27),
(3.28) is equivalent to the formula
$$
\overline{u\circ v}=
\left\{
\begin{array}{ll}
{\bar v}\circ {\bar u} & \forall u,v\in U_{11},\cr
{\bar v}\circ {\bar u}-2{\bar u}\circ {\bar v} &
\forall u\in U_{11},v\in U_{13},\cr
-3{\bar v}\circ {\bar u}+2{\bar u}\circ {\bar v} &
\forall u\in U_{13},v\in U_{11},\cr
{\bar v}\circ {\bar u} &
\forall u,v\in U_{13}.
\end{array}
\right.\eqno(3.30)$$
\par
\vskip 3mm
{\it Proof.}\quad
\rm
We rewrite the equation (3.9) under conditions
(3.26), (3.27).
Then it looks as follows:
$$
Q(u\circ v)=R(Q^{-1}v\circ Qu)+Q(u)\circ {\bar v}-
Q^{-1}(v)\circ RQu+{\bar v}\circ
Qu-{\bar u}\circ Qv.\eqno(3.31)$$
In two cases where
$u,v\in U_{11}$
and $u,v\in U_{13}$
 formula (3.31) becomes very simple:
$$
\overline{u\circ v}={\bar v}\circ {\bar u}.$$
\par
To consider the remaining cases we 
need the following formulas which immediately
follow from (3.28)
$$
R^{-1}(a\circ c)=c\circ a\quad \;
\forall a\in U_{13},\,c\in U_{11},\eqno(3.32)$$
$$
R^{-1}(a\circ c)={4\over 3}a\circ c-{1\over 3}c\circ a\quad
\forall a\in U_{11},\,c\in U_{13}.\eqno(3.33)$$
\par
Now let  $u\in U_{11}$ and $v\in U_{13}$ in (3.31).
We obtain 
$$
Q(u\circ v)=-2{\bar u}\circ {\bar v}+
{2\over 3}{\bar v}\circ {\bar u}+{1\over 3}R({\bar v}\circ {\bar u})$$
or using  formula (3.28)
$$
Q(u\circ v)=2{\bar v}\circ {\bar u}-3{\bar u}\circ  
{\bar v}\ \;\forall u\in U_{11},\;v\in U_{13}.$$
Acting on both sides by $R^{-1}$
and using (3.32), (3.33)
we obtain
$$
\overline{u\circ v}={\bar v}\circ {\bar u}-2{\bar u}\circ {\bar v}\; 
\forall u\in U_{11},\;v\in U_{13}.$$
\par
At last we consider the case
$u\in U_{13}$ and $v\in U_{11}.$
Formula (3.31) gives in this case
$$
Q(u\circ v)=3R({\bar v}\circ {\bar u})+2{\bar u}\circ {\bar v}
-6{\bar v}\circ {\bar u}$$
or using (3.28)
$$
Q(u\circ v)=5{\bar u}\circ {\bar v}-6{\bar v}\circ {\bar u},\ 
\forall u\in U_{13},v\in U_{11}.$$ 
Acting on both sides by
$R^{-1}$
and using (3.32), (3.33) we obtain
$$
\overline{u\circ v}=-3{\bar v}\circ {\bar u}+2{\bar u}\circ {\bar v}
\; \forall u\in U_{13},v\in U_{11}.$$
Thus the formula (3.30) is proved for all cases.
\par
It is evident that starting with
(3.30) we can repeat all considerations backwards for all
cases and come to  formula (3.31).
The lemma is proved.
\par
To formulate the theorem below we recall the following well-known lemma
which we will use without proof.
\par
\vskip 3mm
{\bf Lemma 3.3.}\quad
\it
Let $C(x,y)$ be an algebra defined on
the space $W$ and $*$ be an involutive automorphism
of $C(x,y)$, i.e. $(x^*)^*=x$ and
$$
(C(x,y))^{*}=C(x^{*},y^{*}).\eqno(3.34)$$
Then the space $W$ is a direct sum
$$
W=W_{+}\oplus W_{-}\eqno(3.35)$$
such that
$$
x^{*}=
\left\{\matrix{x\ {\rm if}\ x\in W_{+},\cr
-x\ {\rm if}\ x\in W_{-},\cr}\right.
\eqno(3.36)$$
and
$$
C(W_{+},W_{+})\subset W_{+},\;C(W_{-},W_{-})\subset W_{+},$$
$$
C(W_{+},W_{-})\subset W_{-},\;C(W_{-},W_{+})\subset W_{-}.\eqno(3.37)$$
\par
\vskip 3mm
{\bf Theorem 3.1}\quad
Let $U$ be the space (3.16) and the linear operations
$L,R,Q$ be defined by
(3.15), (3.17), (3.18).
Then the system of equations (3.1)-(3.14) has the following solution:
the 
operator 
$(xey)$
is of the form
$$
(xey)\equiv x\circ y=A_{1}(x,y)+A_{3}(x,y)\ \;\;\forall x,y, \eqno(3.38)$$
where
$$A_{1}(x,y)\in U_{11},\;A_{3}(x,y)\in U_{13}\;\;\forall x,y$$
and
the operator $A_{1}(x,y)$
is symmetric: 
$$ 
A_{1}(x,y)=A_{1}(y,x)\ \;\;\forall x,y,\eqno(3.39)$$
while the operator
$A_{3}(x,y)$ has the property
$$
(-3)^{p(y)}A_{3}(x,y)=-(-3)^{p(x)}A_{3}(y,x)\;\;\forall x,y,\eqno(3.40)$$
i.e.
$$
A_{3}(x,y)=
\left\{
\begin{array}{ll}
-A_{3}(y,x) & \forall x,y\in U_{11},\cr
-A_{3}(y,x) & \forall x,y\in U_{13},\cr
3A_{3}(y,x) & \forall x\in U_{13},
y\in U_{11}.
\end{array}
\right.\eqno(3.41)$$
\par
Moreover, the transformation $\sim$ defined in (3.20)
 is an involutive automorphism of the algebra with the multiplication  $\circ$,
that is
$$
\widetilde{x\circ y}={\tilde x}\circ {\tilde y}\ 
\forall x,y.\eqno(3.42)$$
Hence
$$
U_{+}\circ U_{+}\subset U_{+},\;
U_{-}\circ U_{-}\subset U_{+},\;
U_{-}\circ U_{+}\subset U_{-},\;
U_{+}\circ U_{-}\subset U_{-},
\eqno(3.43)$$
where
$$
U_{+}=U_{11}^{+}+U_{13}^{-},\;
U_{-}=U_{11}^{-}+U_{13}^{+}.$$
The bilinear operators
$(exy)$
and $(xye)$ are 
defined by formulas
(3.26), (3.27).
\par
Conversely,
let $U$ be  a linear space presented in the form (3.16) and three
linear  operators
$L,R,Q$ be
defined by formulas
(3.15), (3.17), (3.18).
\par
Then any bilinear operator 
$x\circ y=A_{1}(x,y)+A_{3}(x,y)$
with $A_{1}(x,y)\in U_{11},A_{3}(x,y)\in U_{13}$
 such that
\par
1)\ transformation (3.20) is an involutive automorphism
of $x\circ y,$
\par
2)\ $A_{1}(x,y),\;A_{3}(x,y)$ satisfy (3.39) and (3.40),
\par
\noindent 
defines a solution of the system of equations
(3.1)-(3.14)
according to formulas
(3.23),
(3.26),
(3.27).
\par
\vskip 3mm
Proof.\quad
\rm
Consider an unknown bilinear operator
$(xey)$ in the form (3.38).
Then  equation (3.11) is equivalent to (3.39).
\par
To prove (3.40) we can use formula (3.28) according to Lemma 3.1.
\par
First we consider the case
$u,v\in U_{11}.$          
Then formula
(3.28) gives
$$
R(u\circ v)=2u\circ v-v\circ u,\ \forall u,v\in U_{11}.$$
Changing places of $v$ and $u$ and adding 
both equalities we have
$$
R(u\circ v+v\circ u)=
u\circ v+v\circ u\ \forall u,v\in U_{11},\eqno(3.44)$$
i.e.,
$u\circ v+v\circ u\in U_{11}.$
The last result is equivalent to (3.40) for $u,v\in U_{11}.$
\par
Consider $u,v\in U_{13}.$
Then again
$$
R(u\circ v)=2u\circ v-v\circ u\ \forall u,v\in U_{13}.$$
In the same way as in the previous case we obtain that (3.40) is true
for $u,v\in U_{13}.$
\par
Now let $u\in U_{13},v\in U_{11}.$
Then considering the element
$u\circ v-3v\circ u,$
we obtain from (3.28) that
$$
R(u\circ v-3v\circ u)=
u\circ v-3v\circ u,\ \forall u\in U_{13},\,v\in U_{11}.\eqno(3.45)$$
This means that
$u\circ v-3v\circ u \in U_{11}.$
This is equivalent to the remaining case in
(3.40).
Thus (3.40) is proved.
\par
We note that doing all considerations above backwards one
can prove that identities (3.39) and  (3.40)
together imply formulas (3.28).
\par
Now we start  proving formula (3.42).
We will use Lemma 3.2 and show that  formulas (3.30)
together are equivalent to the assertion that the transformation
$\sim$ is an involutive automorphism of the bilinear operation
$\circ$.
\par
Consider the case $u,v\in U_{11}.$
According to (3.19), (3.20), (3.21)  formula (3.30) looks as follows:
$$
\widetilde{A_{1}(u,v)}-\widetilde{A_{3}(u,v)}=
A_{1}({\tilde v},{\tilde u})+A_{3}({\tilde v},{\tilde u}),
\ \forall u,v\in U_{11}.$$
Using (3.39), (3.41)
we obtain
$$
\widetilde{A_{1}(u,v)}-\widetilde{A_{3}(u,v)}=
A_{1}({\tilde u},{\tilde v})-A_{3}({\tilde u},{\tilde v}).$$
This is equivalent to 
$$
\widetilde{A_{1}(u,v)}+\widetilde{A_{3}(u,v)}=
A_{1}({\tilde u},{\tilde v})+A_{3}({\tilde u},{\tilde v}),\eqno(3.46)$$
which means
$$
\widetilde{u\circ v}={\tilde u}\circ {\tilde v}\ \forall u,v\in U_{11}.$$
\par
In the same way we obtain from (3.30) the formula
$$
\widetilde{u\circ v}={\tilde u}\circ {\tilde v}\ \forall u,v\in U_{13}.$$
Now we will consider the case $u\in U_{11},v\in U_{13}.$
Formula (3.30) in this case  expressed in terms of 
$A_{1},A_{3}$
looks as follows
$$
\widetilde{A_{1}(u,v)}-\widetilde{A_{3}(u,v)}=
-A_{1}({\tilde v},{\tilde u})-
A_{3}({\tilde v},{\tilde u})
+2A_{1}({\tilde u},{\tilde v})+2A_{3}({\tilde u},{\tilde v})$$
$$
\forall u\in U_{11},v\in U_{13}.$$
Again using (3.39) and (3.41)
we obtain the equation of the form (3.46) which is equivalent to
$$
\widetilde
   {u\circ v}={\tilde u}\circ {\tilde v}\ \;
\forall u\in U_{11},v\in U_{13}.$$
At last we write down (3.30) in the case
$u\in U_{13},v\in U_{11}.$
We obtain
$$
\widetilde{A_{1}(u,v)}-\widetilde{A_{13}(u,v)}=
3A_{1}({\tilde v},{\tilde u})+
3A_{3}({\tilde v},{\tilde u})-
2A_{1}({\tilde u},{\tilde v})-2A_{3}({\tilde u},{\tilde v}),$$
$$
\forall u\in U_{13},v\in U_{11}.$$
Using again (3.39) and (3.41) we come to the relation
$$
\widetilde{u\circ v}={\tilde u}\circ {\tilde v}\ \;
\forall u\in U_{13},v\in U_{11}.$$
Hence we have proved the formula (3.42) for all cases.
Conversely, formula
(3.42) together with formulas
(3.28) imply formulas (3.30) if one repeats the above procedure backwards.
\par
Taking into account Lemma 3.3
we see that relations (3.42) are also proved.
\par
To finish the proof of the first part of the 
theorem we have to prove that all
identities (3.1)-(3.14) are
fulfilled for the considered
solution.
This will be done in the proof
of the converse assertion.
\par
To prove the converse assertion we have
to prove that all formulas
(3.1)-(3.14) are true.
First we note that formulas
(3.1)-(3.3) are identities if
$L$ is the identity operator.
Now we start with formulas
(3.7), (3.10) and
(3.6), (3.8).
\par
Taking into account that in our case $L$
is the identity operator,
one  can  easily see that formulas (3.7) and
(3.10)  coincide with (3.26).
Similarly  formulas (3.6), (3.8)  coincide with (3.27) if 
one takes into account the  
formula (3.26).
\par
According to two 
remarks in the first part of the proof
 we can use formulas (3.28) and
(3.30).
Thus it follows from Lemmas 3.1 and 3.2 that equations (3.9) and
(3.14) hold.
\par
Equation (3.11) is equivalent to relation (3.33).
\par
Equation (3.12) written with the help of (3.26) and (3.27) is
$$
R(Q^{-1}(u)\circ a)+
a\circ {\bar u}-Q^{-1}(u)\circ Ra={\bar u}\circ a\ 
\forall a\in U_{11}.\eqno(3.47)$$
\par
Considering the cases where
$u\in U_{11}$
and $u\in U_{13}$
and changing $Q^{-1}(u)$ by ${\bar u}$ in the case
$u\in U_{11}$
and $Q^{-1}(u)={1\over 3}{\bar u}$
in the case
$u\in U_{13},$
we see that (3.47) is equivalent to two equations
$$
R({\bar u}\circ a)=2{\bar u}\circ a-
a\circ{\bar u} \quad \forall \ a,u\in U_{11},$$
$$
R({\bar u}\circ a)=4{\bar u}\circ a-
3a\circ{\bar u}\quad \forall \ a\in U_{11},u\in U_{13}.$$
\par
But these equations are special sub-cases of formula (3.28).
Thus (3.12) holds.
\par
The formula (3.13) for
$a\in U_{11}$
is a consequence of (3.12).
To prove it in the case
$a\in U_{13},$
we rewrite it using formulas
(3.26) and (3.27):
$$
(R-1)(R(Q^{-1}(v)\circ a)+
a\circ {\bar v}-3Q^{-1}(v)\circ a-
{\bar v}\circ a)-$$
$$
-R({\bar v}\circ a)-a\circ Q(v)+
3{\bar v}\circ a+
Q(v)\circ a)=2{\bar v}\circ a.\eqno(3.48)$$
\par
In the same way as for (3.12) we consider separately the 
cases $v\in U_{11}$
and $v\in U_{13}.$
Changing in the first case
$Q(v)=Q^{-1}(v)$
by ${\bar v}$  and in the second
case $Q(v)$ by
$3{\bar v}$ and $Q^{-1}(v)$
by ${1\over 3}{\bar v}$ and
using formulas (3.28) one comes to the identities.
\par
Formula (3.5) coincides with (3.9) if one denotes 
$x={\bar u}$ and takes into account (3.26).
\par
The last formula we have to prove is (3.4).
Using  formula (3.27) and taking into account that
$L=1$ we rewrite (3.4) as
$$
\overline{(xye)}=(yxe)\eqno(3.49)$$
or
$$
\overline{R(Q^{-1}(y)\circ x)+
x\circ {\bar y}-Q^{-1}(y)\circ Rx}
=$$
$$
=R(Q^{-1}(x)\circ y)+y\circ {\bar x}-Q^{-1}
(x)\circ Ry.\eqno(3.50)$$
\par
Using formulas (3.28) and (3.30),
one can check that formula (3.49) is
valid for all possible four sub-cases
$x,y\in U_{11},U_{13}.$
The theorem is proved.
\par
\vskip 3mm
Now we will consider bilinear
equations obtained from the condition of weak 
commutativity (1.41).
Considering the polarization of equation (1.41)
we obtain
$$
(z,y,(uxx)+(xux)+(xxu))+
(u,y,(zxx)+(xzx)+(xxz))+
(x,y,\{(xzu)\}=$$
$$
=((u,x,x)+(xux)+(xxu),y,z)+((zxx)+(xzx)+
(xxz),y,u)+(\{(xzu)\},y,x),\eqno(3.51)$$
where the symbol
$\{(xzu)\}$
means the symmetrization of
$(xzu),$ i.e.
$$
\{(xzu)\}=
(xzu)+(zxu)+(xuz)+
(zux)+(uxz)+(uzx).\eqno(3.52)$$
Substituting in (3.51) first
$x=u=e$
and then
$x=y=e$
we obtain two equations
$$
(z,y,e)+(e,y,(R+Q+1)z)=
(eyz)+((R+Q+1)zye),\eqno(3.53)$$
$$
(ze(R+Q+1)u)+
(ue(R+Q+1)z)+
\{(euz)\}=$$
$$
=((R+Q+1)uez)+((R+Q+1)zeu)+
R\{(euz)\}.\eqno(3.54)$$
\par
We recall that it follows from Theorem 1.2
that in the  weakly commutative case
the subspace
$U_{13}^{+}=0.$
 So in this case instead of the
space (3.16) one has to consider  the space
$$
U=U_{11}^{+}\oplus U_{11}^{-}+U_{13}^{-}.\eqno(3.55)$$
It turns out that this condition is not only necessary but
also sufficient.
\par
\vskip 3mm
{\bf Theorem 3.2.}\quad
\it
Let $U$ be the space
(3.55) and  the linear operations $L,R,Q$
be defined by (3.15), (3.17), (3.18).
Then the system of equations
(3.1)-(3.14),
(3.53), (3.54) has the same solutions
as in Theorem 3.1 
defined by formulas
(3.38)-(3.40)
and (3.23), (3.26), (3.27).
\par
Conversely,
let  $U$ be a linear space given by (3.55)
with operators $L,R,Q$
defined by formulas (3.15), (3.17), (3.18).
Then any bilinear operator
$x\circ y=A_{1}(x,y)+A_{3}(x,y)$
with 
$A_{1}(x,y)\in U_{11},\;A_{3}(x,y)\in U_{13}$
such that
\par
1)\ the transformation (3.20) is an involutive automorphism,
\par
2)\ $A_{1}(x,y),A_{3}(x,y)$
satisfy (3.39) and (3.40),
\par
\noindent
defines a solution of the system of 
equations (3.1)-(3.14),
     (3.53), (3.54) by formulas (3.23), (3.26), (3.27).
\par
\vskip 3mm
Proof.\quad
\rm
Due to Theorem 3.1
we only  have to prove that equations
(3.53), (3.54) are consequences of (3.1)-(3.14) 
under the condition
$U_{13}^{+}=0.$
\par
Let us rewrite (3.53) as
$$
(ey(R+Q)z)=((R+Q)zye).\eqno(3.56)$$
We note that in the case $U_{13}^{+}=0$ 
the vector $(R+Q)z\in U_{11}.$
Thus the equation (3.56) is equivalent to
$$
(eyz)=(zye)\quad \forall z\in U_{11}.\eqno(3.57)$$
But the last equation coincides with (3.12).
\par
Now we rewrite (3.54) as
$$
(R-1)\{(euz)\}=
z\circ(R+Q)u-
(R+Q)u\circ z+u\circ(R+Q)z-
(R+Q)z\circ u.\eqno(3.58)$$
It follows from (3.17) that the left-hand side
belongs to $U_{13}.$
\par
The same it true for the right-hand side according to the
symmetry on $u,z$
and the formula (3.39).
Thus one has to consider only the part $A_{3}(x,y)$
 of the multiplication 
$x\circ y.$
\par
Using the formulas (3.26), (3.27) and then (3.29)
we rewrite the left-hand side of (3.58) as follows
$$
(R-1)\{(euz)\}=(R-1)({\bar u}\circ z+
{\bar z}\circ u+u\circ z+z\circ u
+R(Q^{-1}(u)\circ z)+$$
$$
z\circ {\bar u}-Q^{-1}(u)\circ Rz+
+R(Q^{-1}(z)\circ u)+u\circ {\bar z}-Q^{-1}(z)\circ Ru)=
$$
$$
=2{\bar u}\circ z+2{\bar z}\circ u+u\circ z+z\circ u+
(R-1)z\circ Q^{-1}u-2Q^{-1}u\circ (R-1)z
$$
$$
+(R-1)u\circ Q^{-1}z-2Q^{-1}z\circ(R-1)u.\eqno(3.59)$$
\par
 Let $x=x_{1}+x_{3}$ where $x_{1}\in U_{11},\ x_{3}\in U_{13}.$
Then the right-hand side of (3.58) can be rewritten as follows
$$
z\circ(u_{1}+\overline{u_{1}})-(u_{1}+\overline{u_{1}})
\circ z+u\circ (z_{1}+\overline{z_{1}})-(z_{1}+\overline{z_{1}})
\circ u,\eqno(3.60)$$
because
$(R+Q)z_{3}=0,$
when $U_{13}^{+}=0.$
\par
It follows from (3.60) that the right-hand side equals zero when 
$u,z\in U_{13}.$
But the left-hand side is also equal to zero in this case.
Indeed,
$\overline{z_{3}}=-z_{3},
\overline{u_{3}}=-u_{3}$
and the left-hand side becomes symmetric on
$u_{3},z_{3}$
while 
$A_{3}(z_{3},u_{3})$
is a skew-symmetric operator
(see (3.41)).
\par
In the same way we consider the case
$u,z\in U_{11}.$
The operator
$A_{3}(u_{1},z_{1})$
is also skew-symmetric.
Hence the $U_{13}$ component
of (3.60) equals  the 
$U_{13}$ component of
$$
{1\over 2}[(z_{1}-\overline{z_{1}})\circ (u_{1}+\overline{u_{1}})-
(u_{1}+\overline{u_{1}})\circ (z_{1}-\overline{z_{1}})+$$
$$
(u_{1}-\overline{u_{1}})\circ (z_{1}+\overline{z_{1}})-
(z_{1}+\overline{z_{1}})\circ (u_{1}-\overline{u_{1}})].\eqno(3.61)$$
But according to (3.43) the
$U_{13}$ component of  the products
$u_{11}^{+}\circ z_{11}^{-}$
and $u_{11}^{-}\circ z_{11}^{+}$
is zero.
Thus (3.61) is zero.
\par
Similar considerations show that  the
$U_{13}$ component of (3.59)
is also zero for 
$u,z\in U_{11}.$
\par
For the remaining cases it is enough to consider
$u\in U_{11},\ z\in U_{13}$ because both sides are
symmetric with respect to  $u$ and $z$.
Also we can take $u_{1}\in U_{11}^{+}$ because 
$A_{3}(U_{11},U_{13})=0$
when $U_{13}^{+}=0$
(see (3.43)).
Under these conditions the $U_{13}$ component of (3.59)
equals  the $U_{13}$
component of
$$
2(2u_{1}\circ z_{3}-
2z_{3}\circ u_{1}+u_{1}\circ z_{3}+z_{3}\circ u_{1}+
2z_{3}\circ u_{1}-4u_{1}\circ z_{3})=$$
$$
=2(-u_{1}\circ z_{3}+z_{3}\circ u_{1})=4u_{1}\circ z_{3}.$$
\par
In the last equality we have used  formula (3.41).
By reasoning in the same way the right-hand side also equals
$$
2z_{3}\circ u_{1}-2u_{1}\circ z_{3}=4u_{1}\circ z_{3}.$$
\par
Thus the equality (3.54) is true and the theorem is proved.
\par
The following theorem is a corollary of Theorem 3.1 and
Theorem 3.2.
\par
\vskip 3mm
{\bf Theorem 3.3.}\quad
\it
Let a tripotent $e$ of a generalized Jordan triple system of 
the second order be a  
left unit.
Then the space $U$ of the triple system is a direct sum
$$
U=U_{11}^{+}\oplus U_{11}^{-}\oplus 
U_{13}^{+}\oplus U_{13}^{-},\eqno(3.62)$$
and the triple system has the form
$$
(xyz)=(Q^{-1}(y)\circ x)\circ z+x\circ (Q^{-1}(y)\circ z)
-Q^{-1}(y)\circ (x\circ z),\eqno(3.63)$$
where $Q(x)$ is defined by (3.18).
\par
The bilinear operator $x\circ y$ has the form
$$
x\circ y=A_{1}(x,y)+A_{3}(x,y),\eqno(3.64)$$
where
$A_{1}(x,y)\in U_{11},\ A_{3}(x,y)\in U_{13},$
and have the following properties
$$
A_{1}(x,y)=A_{1}(y,x)\quad \forall x,y,\eqno(3.65)$$
$$
A_{3}(x,y)=
\left\{
\begin{array}{ll}
-A_{3}(y,x)& \forall x,y\in U_{11},\cr
-A_{3}(y,x)& \forall x,y\in U_{13},\cr
3A_{3}(y,x)& \forall x\in U_{13},
y\in U_{11}.
\end{array}
\right.\eqno(3.66)$$
\par
Moreover, the transformation
$$
{\tilde x}=
\left\{
\begin{array}{cl}
x & \forall x\in U_{+}=U_{11}^{+}\oplus U_{13}^{-},\cr
-x & \forall x\in U_{-}=U_{11}^{-}\oplus U_{13}^{+}.
\end{array}
\right.\eqno(3.67)$$
is an involutive automorphism of $x\circ y$.
Hence 
$$
U_{+}\circ U_{+}\subset U_{+},\;
U_{-}\circ U_{-}\subset U_{+},\;
U_{-}\circ U_{+}\subset U_{-},\;
U_{+}\circ U_{-}\subset U_{-}.
\eqno(3.68)$$
\par 
If in addition the triple system
is weakly commutative,
then the 
subspace
$U_{13}^{+}$ is equl to zero.
Thus
$$
U=U_{11}^{+}\oplus U_{11}^{-}\oplus U_{13}^{-},\eqno(3.69)$$
while the multiplication
$x\circ y$ has 
 the same properties 
 (3.65)-(3.68) as in the general case.}
\par
\rm
In \S 2 we gave examples of the Peirce decomposition
defined by a left unit.
The triple system 
$A_{kn}-A_{nk}$
illustrates the general case and the
triple system $D_{nk}$ illustrates the
weakly commutative case
(see (2.13)).
\newpage
\begin{center}
\bf References
\end{center}
\vskip 3mm
\rm
\begin{enumerate}
\item
B. N. Allison,
A class of nonassociative algebras with
involution containing the class of
Jordan algebras, Math.Ann.237(1978), 133-156,
\item
B. N. Allison,
Appendix, Peirce decompositions in structurable
algebras, in the paper: B.Allison, G.Benkart, Y.Gao, Lie algebras
graded by the root systems $BC_r,\;r \ge 2$ , to appear,
\item
B. N. Allison, J.R. Faulkner,
Elementary Groups And Invertibiolity for Kantor Pairs, Comm. Alg.,
27(2),(1999),519-516,
\item
N. Jacobson,
Lie and Jordan triple systems,
Amer.J.Math.71(1949), 149-170,
\item
 J.R. Faulkner,  Structuruble triples, Lie triples, and symmetric
 spaces,
 Forum Math., 6,(19949, 637-650,
\item
N. Kamiya, A structure theory of Freudenthal-Kantor
triple systems,
J.Alg.
110
(1987)
108
-
123,
\item
N. Kamiya, A structure theory of Freudenthal-Kantor triple 
systems
II.Comm.
Math.
Univ.
Sancti 
Pauli 
38
(1989)
41
-
60,
\item
N. Kamiya,
A structure theory of Freudenthal-Kantor triple systems III.
Mem.Fac.
Sci.
Shimane Univ.23(1989), 33-51,
\item
N. Kamiya,
On Freudenthal-Kantor triple systems and 
generalized structurable algebras,
Proceeding of the International
conference of
nonassociative algebras
and its applications,
Mathematics and
its
 Applications
303
(1994)
198-203,
Kluwer Academic
Publisher,
\item
N. Kamiya, On the Peirce decompositions for Freudenthal-Kantor
triple systems,
Comm.
Alg.,
25
(6)
(1997), 
1833-1844,
\item
I. L. Kantor, 
Models of exceptional Lie algebras, Soviet Math.Dokl.14
(1973), 254-258,
\item
I. L. Kantor,
Some generalizations of Jordan 
algebras,
Trudy 
Sem.
Vektor.
Tensor.
Anal.
16
(1972),
407-499
(Russian),
\item
I. L. Kantor, A generalization of the 
Jordan 
approach 
to
symmetric
Riemennion 
spaces,
The 
Monster
and 
Lie 
algebras, 
Ohio 
University,
Mathematical 
Resarch 
Institute 
Publications 7,
(1998), 221-234,
\item
O. Loos,
Jordan Pairs,
Lecture Note in 
Math.
406(1975)
Springer.
\end{enumerate}
\end{document}